\documentclass[11pt,a4paper,dvipsnames,eqrno]{article}
\usepackage[utf8]{inputenc}
\usepackage{a4wide}

\usepackage{graphicx}
\usepackage{amssymb}
\usepackage{amsmath}
\usepackage{amsthm}
\usepackage{thmtools}
\usepackage{mathtools}
\usepackage{multirow}
\usepackage{blkarray}
\usepackage[dvipsnames]{xcolor}
\usepackage{etoolbox}
\usepackage{enumerate}
\usepackage[shortlabels]{enumitem}
\usepackage{verbatim}
\usepackage{comment}
\usepackage{array}
\usepackage{dsfont}
\usepackage{soul}
\usepackage{braket}
\usepackage{standalone}
\usepackage[ruled,linesnumbered]{algorithm2e}

\usepackage{relsize}
\usepackage[pdftex, colorlinks, linkcolor=blue, citecolor=blue, urlcolor=blue]{hyperref}
\usepackage[nameinlink]{cleveref}
\hypersetup{linkcolor=MidnightBlue, citecolor=PineGreen, urlcolor=Emerald}
\usepackage{pgfplots}
\usepackage[margin=1cm]{caption}
\usepackage{subcaption}
\pgfplotsset{width=7cm, compat=1.10}
\usepackage{todonotes}

\newcommand{\A}{\mathcal{A}}
\newcommand{\B}{\mathcal{B}}
\newcommand{\R}{\mathbb{R}}
\newcommand{\Z}{\mathbb{Z}}
\newcommand{\C}{\mathbb{C}}

\newcommand{\CP}{\mathbb{CP}}

\newcommand{\Q}{\mathbb{Q}}

\newcommand{\Ical}{\mathcal{I}}
\newcommand{\Jcal}{\mathcal{J}}
\newcommand{\Acal}{\mathcal{A}}
\DeclareMathOperator{\res}{res}
\DeclareMathOperator{\Sing}{Sing}
\DeclareMathOperator{\Proj}{Proj}

\newcommand{\csm}{c_{SM}}
\theoremstyle{definition}
\newtheorem{theorem}{Theorem}[section]
\newtheorem{corollary}[theorem]{Corollary}
\newtheorem{proposition}[theorem]{Proposition}
\newtheorem{definition}[theorem]{Definition}

\newtheorem{example}[theorem]{Example}
\newtheorem{lemma}[theorem]{Lemma}
\newtheorem{remark}[theorem]{Remark}

\usepackage{graphicx} % Required for inserting images

\title{Hypersurface Arrangements with Generic Hypersurfaces Added}
\author{Bernhard Reinke, Kexin Wang}
\date{}

\begin{document}
\maketitle
\abstract{
    The Euler characteristic of a very affine variety encodes the number of critical points of the likelihood equation on this variety.
    In this paper, we study the Euler characteristic of the complement of a hypersurface arrangement with generic hypersurfaces added.
    For hyperplane arrangements, it depends on the characteristic polynomial coefficients and generic hypersurface degrees. 
    As a corollary, we show that adding a degree-two hypersurface to a real hyperplane arrangement enables efficient sampling of a single interior point from each region in the complement.
    We compare the method to existing alternatives and demonstrate its efficiency.
    For hypersurface arrangements, the Euler characteristic is expressed in terms of Milnor numbers and generic hypersurface degrees.
    This formulation further yields a novel upper bound on the number of regions in the complement of a hypersurface arrangement.}

\section{Introduction}
In this article, we study the topological Euler characteristic
\begin{equation}
\chi \bigl(\C^n\setminus(\Acal\cup V(\prod_{i=1}^r{g_i}) )\bigr) 
\label{eq:chi}
\end{equation}
where $\Acal = \{V(f_1),\ldots,V(f_k)\}$ is a fixed hypersurface arrangement in $\C^n$ with $f_i \in \C[x_1,\ldots,x_n]$ for $i=1,\ldots,k$ and $g_i\in \C[x_1,\ldots,x_n]$ is a generic polynomial of degree $d_i$ for $i=1,\ldots,r$.
When $\C^n\setminus \Acal$ is smooth and very affine, by \cite{huh2013maximum}, $(-1)^n\chi \bigl(\C^n\setminus(\Acal\cup V(\prod_{i=1}^r{g_i})\bigr)$ coincides with the number of critical points of
\begin{equation}
\psi(x_1,\ldots,x_n) = \sum_{i=1}^k u_i \log f_i(x) + \sum_{i=1}^r v_i \log g_i(x),
\label{eq:loglikelihood}
\end{equation}
for generic parameters $u=(u_1,\ldots,u_k)\in \C^k$ and $v=(v_1,\ldots,v_d)\in \C^d$.
In particular, when $V(f_1),\ldots, V(f_k)$ are hyperplanes, $\C^n\setminus \Acal$ is smooth and very affine when the dimension of the space spanned by the normals of $f_1,\ldots,f_k$ has dimension $n$ (in other words, the hyperplane arrangement $\Acal$ is essential).
The function $\psi$ arises as the likelihood function in algebraic statistics:

Consider the $(k+r-1)$-state discrete distributions $(f_1(x),\ldots,f_k(x),g_1(x),\ldots,g_r(x))$ parameterized by $x\in \R^n$ and the observed count data $(u_1,\ldots,u_k,v_1,\ldots,v_r)$. 
The problem of computing the maximum likelihood estimation (MLE) corresponds to maximizing the function 
$\psi(x_1, \ldots, x_n).$
In particular, the MLE is one of the critical points for $\psi(x)$. 
The total number of complex critical points of $\psi(x)$, known as the ML degree \cite{catanese2006maximum}, quantifies the algebraic complexity of determining the maximum likelihood estimation.

%\KW{Connection to physics:} 
%In particle physics, such function $\psi(x)$ is the scattering potential, with one notable example coming from the Cachazo–He–Yuan (CHY) model \cite{CHY}. 

%When $\mathcal{A}$ is a hyperplane arrangement, the characteristic polynomial $\chi_\mathcal{A}(t)$ is a polynomial that encodes the combinatorial structure of the arrangement. 
%It is computed based on the intersection lattice of the hyperplanes. 
%When the hyperplanes in $\mathcal{A}$ are defined by polynomials with real coefficients, we have:
%$$|\chi_\mathcal{A}(1)| = \text{(number of bounded regions in } \mathbb{R}^n \setminus \mathcal{A}),$$
%and 
%$$|\chi_\mathcal{A}(-1)| = \text{(number of regions in } \mathbb{R}^n \setminus \mathcal{A}).$$

In this paper, we first show that when $\A$ is a hyperplane arrangement, the Euler characteristic \eqref{eq:chi} 
only depends on the coefficients of the characteristic polynomial $\chi_\mathcal{A}(t)$ and the degrees $d_1, \ldots, d_r$ of the generic hypersurfaces.

\begin{theorem}\label{thm: hyperplane multiple hypersurface}
    Let $\A = \{V(f_1),\ldots,V(f_k)\} \subset \C^n$ be a hyperplane arrangement with characteristic polynomial $\chi_A(t) = \sum^n_{i=0} a_i t^i$. 
    Let $d_1, \dots d_r \in \Z_+$.
    Denote $b_i[d_1, \dots, d_r]$ to be the coefficient of $z^i$ in  $$\frac{1}{(1-z)} \prod^r_{i = 1} \frac{1- z}{1 + (d_i - 1) z}.$$
    Then, for generic forms $g_1, \dots, g_r$ on $\C^n$ of degree $d_1, \dots, d_r$, we have: 
    $$\chi(\C^n \setminus (\A  \cup V(\prod_{i=1}^rg_i))) = \sum^n_{i=0} a_i b_i[d_1, \dots, d_r].$$
\end{theorem}

When a single generic hypersurface $V(g)$ of degree $d$ is added, the above theorem implies that the Euler characteristic \eqref{eq:chi} equals $\chi_\Acal(1-d)$.

For $d= 0$, $g$ is a constant function. 
We recover the classical fact that $\chi(\C^n \setminus \A) = \chi_\Acal(1)$. When $\Acal$ is a real and essential hyperplane arrangement, we also have by \cite{zaslavsky1975facing} that $ \chi_\Acal(1)$ counts the number of bounded regions of $\R^n \setminus \Acal$. Varchenko in \cite{varchenko1995critical} relates this to number of complex critical points of $\psi$: if we choose each $u_i>0$, $\psi$ has exactly one critical point in each bounded region.

%This result is proved when $\Acal$ is a real hyperplane arrangement by Varchenko in \cite{varchenko1995critical} and in general by Orlik and Terao in \cite{orlik1995number}.
%When $\Acal$ is real and essential, $\chi(\C^n \setminus \A)$ is the number of complex critical points of $\psi(x) = \sum_{i=1}^k u_i \log f_i(x)$ with generic $u_1,\ldots,u_k$ and $|\chi_\Acal(1)|$ equals the number of bounded regions for $\R^n\setminus \Acal$. 
%If we choose each $u_i>0$, $\psi$ tends to negative infinity at the boundary of each bounded region in $\R^n\setminus\Acal$ ensuring at least one local maximum of $\psi$ in each bounded region.
%However, the number of complex critical points of $\psi$ equals the number of bounded regions so all complex critical points of $\psi$ are real and there is exactly one real critical point per bounded region of $\R^n\setminus\Acal$ \cite[Proposition 11]{catanese2006maximum}.

For $d=2$, the result confirms a numerical observation made by the second author in \cite{breiding2024computing} stated below. 

\begin{theorem}\label{thm: one real point per region}
Suppose $\Acal = \{V(f_1),\ldots,V(f_k)\}$ is an essential real hyperplane arrangement.
% , i.e. there exists some $S\subseteq [k]$ with $\bigcap_{i\in S} V(f_i)$ is a point
Let $g$ be a generic degree two polynomial satisfying $g(x)>0$ for all $x\in \R^n$. 
Then for the function
$$
\psi(x) = \sum_{i=1}^k u_i \log |f_i(x)| -v \log|g(x)|,
$$
where $u_1,\ldots,u_k,v\in \R_{>0}$ and $2v>\sum_{i=1}^k u_i$, all of its critical points are real and there is precisely one critical point per region for $\R^n \setminus \Acal$.
\end{theorem}

The result directly yields an efficient algorithm for sampling interior points of all regions, including unbounded ones, in the complement of a given hyperplane arrangement. 
% Additionally, one can determine which regions are bounded by solving critical points of $\sum_{i=1}^k u_i \log |f_i(x)|$ with $u_1,\ldots,u_k \in \R_{>0}$.
Sampling interior points from each region has applications in optimization, machine learning, automatic control, and other fields.

Sampling from unbounded regions is crucial for applications such as utility maximization and resource allocation, where decision variables, such as prices or investments, may be unbounded. 
Additionally, dual optimization problems often involve unbounded regions, even when the primal problem is bounded. 
This framework also facilitates the study of central hyperplane arrangements, where all regions are unbounded, e.g. \cite{baldi2019polynomial,brandenburg2024best}.

% Each region of $\R^n \setminus \Acal$ corresponds to a sign pattern $\sigma\in \{-1,+1\}^k$, determined by evaluating $(f_1,\ldots,f_k)$ at any interior point of the region. 
% Understanding which sign patterns are feasible, identifying whether regions are bounded or unbounded, and sampling interior points from each region have applications in optimization, machine learning, automatic control and other fields.

The interior point method is a foundational algorithm for solving linear and nonlinear convex optimization problems \cite{nesterov1994interior}. 
The algorithm relies on iteratively moving through the interior of the feasible region along the central path to find the optimal solution, so finding an initial interior point is crucial. 
In practice, one may need to optimize a certain function defined differently over convex regions of some hyperplane arrangement complement \cite{brandenburg2024best,montufar2014number}, where sampling one point from each region is essential.

For a fixed set $S$ of points in $\R^{n-1}$, the separability arrangement \cite{baldi2019polynomial, brysiewicz2023computing} in $\R^n$ encodes the possible ways to partition the set of points into two classes using a hyperplane. 
This concept is central in machine learning algorithms like Support Vector Machines (SVM) and the Perceptron, both of which rely on linear separability.
A special case occurs when $S$ is the set of vertices of the hypercube, known as the threshold arrangement, which has applications in neural networks \cite{zunic2004encoding,cueto2010geometry,montufar2015does}.
Regions in the complement of the separability arrangement are in bijection with linearly separable partitions of $S$ and an interior point from each region represents a hyperplane that achieves the partition.
Another notable family of arrangements is the resonance arrangement, whose regions have practical meanings in the study of economics and quantum field theory. For an overview, see \cite[Section 1]{kuhne2023universality}.

If the hyperplanes in Theorem \ref{thm: one real point per region} are replaced by hypersurfaces with real coefficients, and the condition $2v > \sum_{i=1}^k u_i \deg(f_i)$ is imposed, the guarantee that all complex critical points are real and that there is exactly one critical point per region no longer holds. 
Nevertheless, there is at least one critical point in each region, and the total number of critical points for
$$
\psi(x) = \sum_{i=1}^k u_i \log |f_i(x)| -v \log|g(x)|
$$ 
is an upper bound for the number of regions in $\R^n\setminus\Acal$. 
This motivates a generalization of Theorem \ref{thm: hyperplane multiple hypersurface} from hyperplane arrangement to hypersurface arrangement. 
The Euler characteristic \ref{eq:chi} satisfies a similar formula to that in Theorem \ref{thm: hyperplane multiple hypersurface} with the coefficients of the Euler characteristic replaced by Milnor numbers \cite{milnor1968singular,teissier1973cycles}.

\begin{theorem}\label{thm: hypersurface arrangement with multiple hypersurface}
Let $\A= V(f) \subset \C^n$ be a hypersurface arrangement where $f = \prod_{i=1}^k f_i$. 
We identify $\C^n$ with the affine chart $\CP^n\setminus V(x_0).$
Let $d_1, \dots, d_r \in \Z_+$.
Then for generic forms $g_1, \dots, g_r$ of degree $d_1, \dots, d_r$, we have
$$\chi(\C^n \setminus (\A  \cup V(g_1 \cdots g_r))) = \sum^n_{i=0} (-1)^{n-i} \mu^{n-i}({}^hfx_0) b_i[d_1, \dots, d_r],$$
where ${}^hf$ is the homogenization of $f$ using $x_0$, $\mu^i(\cdot)$ denotes the $i$-th Milnor number of a hypersurface in $\CP^n$ and $b_i[d_1, \dots, d_r]$ is as defined in Theorem \ref{thm: hyperplane multiple hypersurface}.
\end{theorem}

Theorem \ref{thm: hypersurface arrangement with multiple hypersurface} results in the following upper bounds on the number of regions. 
% To the best of the authors' knowledge, this is the first upper bound for the number of regions for hypersurface arrangements in $\R^n$ with general $n$. 

\begin{theorem}\label{thm: number of regions bounded by complex euler}
Suppose $f_1,\ldots,f_k\in \R[x_1,\ldots,x_n]$ and $\C^n\setminus(\Acal \cup V(g))$ is smooth and very affine where $\Acal = \cup_{i=1}^k V(f_i)$ and $g(x_1,\ldots,x_n)\in \R[x_1,\ldots,x_n]$ is a generic degree two polynomial such that $g(x)>0$ for any $x\in \R^n$. Let $f=\prod_{i=1}^k f_i$ and $d=\deg f$.
Then, 
$$
\# \text{ regions of } \R^n \setminus \Acal \leq  |\chi(\C^n - (\Acal\cup V(g))| = \sum_{i=0}^n\mu^i({}^hfx_0) \,
$$
where ${}^hf$ is the homogenization of $f$. 
Moreover, $\mu^i({}^hfx_0) \leq d^i$ for all $i=0,\ldots,n$. Hence,
$$
\# \text{ regions of } \R^n \setminus \Acal \leq \frac{(d+1)^{n+1}-1}{d} = \mathcal{O}(d^n).
$$
\end{theorem}

The second inequality in Theorem \ref{thm: number of regions bounded by complex euler} matches the asymptotic bound established in \cite[Theorem 1]{basu1996bounding}. 
See also \cite{barone2012bounding} for a more generalized result on the number of regions for hypersurface arrangements on varieties.
In contrast, the first inequality in Theorem \ref{thm: number of regions bounded by complex euler} may provide an asymptotically tighter bound than the second, as demonstrated in Example \ref{ex:concentric circles}.

The rest of the paper is organized as follows. 
In Section \ref{sec: hyperplane}, we briefly review the classical results about characteristic polynomials and the Euler characteristic of complete intersections. 
We then prove Theorem \ref{thm: hyperplane multiple hypersurface} and Theorem \ref{thm: one real point per region}. 
We provide the algorithm for computing one interior point per region of the complement of real hyperplane arrangements and we illustrate its speed compared with the other method.
In Section \ref{sec: hypersurface}, we briefly recall Milnor numbers, CSM class.
We prove Theorem \ref{thm: hypersurface arrangement with multiple hypersurface} and then present results about bounding the total number of regions.

\section{Hyperplane arrangement with generic hypersurfaces added}\label{sec: hyperplane}

\subsection{Preliminaries}
This subsection briefly introduces classical results about hyperplane arrangements and the Euler characteristic of complete intersections.

Let $\Acal$ be a hyperplane arrangement in $\mathbb{K}^n$ ($\mathbb{K}$ is a field) and let $L(\Acal)$ be the set of all nonempty intersections of subsets of hyperplanes in $\A$ including $\C^n$ itself as the empty intersection. We order elements in $L(\Acal)$ by reverse inclusion and call $L(\Acal)$ \emph{the intersection poset} of $\Acal$.

For any poset $P$, we denote $\text{Int}(P)$ the set of all closed intervals of $P$. A fundamental invariant of posets with all closed intervals finite is its Möbius function.

\begin{definition}
     The \emph{Möbius function} of a poset $(P,\leq)$ is a function $\mu \colon \text{Int}(P) \rightarrow~\mathbb{Z}$ such that $\mu(x,x)=1$ and $\sum_{z \in [x,y]} \mu(x,z)=0$ for every $x,y \in P$. Note that we write $\mu(x,y):=\mu([x,y])$ and $\mu(x):= \mu([\hat{0},x])$, where $\hat{0}$ is the minimal element of $P$.
 \end{definition}

 We now define the characteristic polynomial of a hyperplane arrangement.

 \begin{definition}
The characteristic polynomial $\chi_\A(t)$ of an hyperplane arrangement $\A$ is 
$$
\chi_\A(t) = \sum_{x\in L(\A)} \mu(x)t^{\dim x}.
$$
 \end{definition}

\begin{example}
\label{ex: hyperplane arrangement}
Consider the real hyperplane arrangement $\Acal\subseteq \R^2$ shown in Figure \ref{fig:arrangement_and_poset} (left). Its intersection poset $L(\Acal)$ is shown in Figure \ref{fig:arrangement_and_poset} (right).

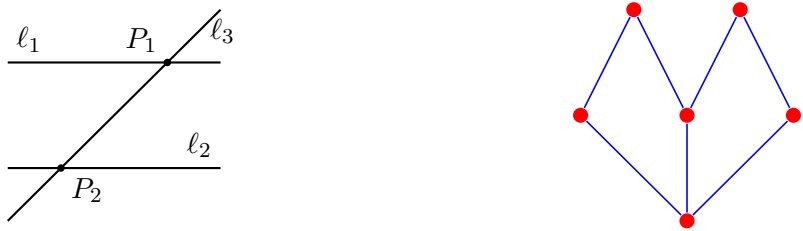
\begin{figure}[htbp]
    \centering
    \begin{minipage}{0.45\textwidth}
        \centering
        \begin{tikzpicture}[scale = 0.7]

  % Draw the first parallel line L1
  \draw[thick] (-2, 2) -- (2, 2) node[pos=0.1, above] {$\ell_1$};

  % Draw the second parallel line L2
  \draw[thick] (-2, 0) -- (2, 0) node[pos=0.9 , above] {$\ell_2$};

  % Draw the third intersecting line L3
  \draw[thick] (-2, -1) -- (2, 3) node[pos=0.9, right] {$\ell_3$};

  % Mark the points of intersection
  \fill (1,2) circle (2pt) node[above left] {$P_1$}; % Intersection of L1 and L3
  \fill (-1,0) circle (2pt) node[below right] {$P_2$}; % Intersection of L2 and L3

\end{tikzpicture}
    \end{minipage}%
    \quad % Horizontal spacing
    \begin{minipage}{0.45\textwidth}
        \centering
        \begin{tikzpicture}[scale=0.7, every node/.style={circle, fill=red, inner sep=2pt}, line width=0.6pt]

\node (P1) at (-1,2) {};
\node (P2) at (1,2) {};
\node (l1) at (-2,0) {};
\node (l3) at (0,0) {};
\node (l2) at (2,0) {};
\node (R2) at (0,-2) {};

\draw[blue] (P1) -- (l1) -- (R2);
\draw[blue] (P2) -- (l2) -- (R2);
\draw[blue] (P1) -- (l3) -- (R2);
\draw[blue] (P2) -- (l3);
\end{tikzpicture}
    \end{minipage}
    \caption{Hyperplane arrangement and its corresponding poset}
    \label{fig:arrangement_and_poset}
\end{figure}

By definition, the Möbius function $\mu$ takes values $\mu(\R^2) = 1, \mu(\ell_1) = \mu(\ell_2) = \mu(\ell_3) = -1, \mu(P_1) = \mu(P_2) = 1$. Hence, the characteristic polynomial for $\Acal$ is $\chi_\Acal(t) = t^2 - 3t +2$.
\end{example}

Thomas Zaslavsky in \cite{zaslavsky1975facing} showed that evaluations of the characteristic polynomial give region count.

\begin{proposition}

Let $\Acal\subseteq \R^n$ be an essential hyperplane arrangement and $\chi_\A(t)$ its characteristic polynomial. Then,
\begin{align*}
\text{(the number of bounded regions in } \mathbb{R}^n \setminus \mathcal{A}) & = (-1)^{\text{rank}\A}\chi_\mathcal{A}(1) \\
\text{(the number of regions in } \mathbb{R}^n \setminus \mathcal{A}) & = (-1)^n \chi_\Acal(-1).
\end{align*}
\end{proposition}

An alternative formula for the characteristic polynomial of an arrangement is due to Hassler Whitney, which we will use in proving Theorem \ref{thm: hyperplane multiple hypersurface}.

\begin{proposition}[Whitney's Theorem](see for example \cite[Theorem 2.4]{StanleyHyperplanes})\label{thm:whitney}

    Let $\A$ be an arrangement in an $n$-dimensional vector space. Then,
\begin{equation*}
\sum_{\substack{\B \subset \A\\ \bigcap \B \text{ non-empty}}} (-1)^{\# \B} z^{\dim \bigcap \B} = \chi_A(z)
    \label{eqn:whitney}
\end{equation*}
\end{proposition}

\begin{example}
We continue with the hyperplane arrangement in Example \ref{ex: hyperplane arrangement}.
We list all non-empty intersections of subsets of $\A$ in Table \ref{table: whitney}. 
From Proposition \ref{thm:whitney}, we obtain 
$\chi_\A(t) = (-1)^0 t^2 - 3 t^1 + 2 (-1)^2 t^0 = t^2 - 3t +2,$ which coincides with the one computed directly from the Möbius function.

\begin{table}[htbp]
\centering
\begin{tabular}{|c|c|c|}
\hline
    $\B$ & $\# \B$ & $\dim \cap \B$ \\
    \hline
    $\emptyset$ & 0 & 2\\
    \hline
    $\{\ell_1\}$ & 1 & 1\\
    \hline
    $\{\ell_2\}$ & 1 & 1\\
    \hline
    $\{\ell_3\}$ & 1 & 1\\
    \hline
    $\{\ell_1,\ell_3\}$ & 2 & 0\\
    \hline
    $\{\ell_2,\ell_3\}$ & 2 & 0\\
    \hline  
\end{tabular}
\caption{non-empty intersections of subsets of $\A$}
\label{table: whitney}
\end{table}

\end{example}

We recall a classical formula for the Euler characteristic of a complete intersection (in this form due to \cite{HirzebruchICM}).
\begin{lemma}\label{lem: hirzeburch}
Let $d_1, \dots, d_r \geq 1$. For generic $g_1, \dots g_r$ forms of degree $d_1, \dots, d_r$ on $\CP^n$, $V(g_1, \dots, g_r)$ is a non-singular complete intersection (or empty if $r > n$). 
Then the Euler characteristic of $V(g_1, \dots, g_r)$ only depends on $d_1, \dots, d_r, n$ and is denoted by $c_n[d_1, \dots, d_r]$
and we have
\begin{equation}
    \sum_{n \geq 0} c_n[d_1, \dots, d_r] z^n = \frac{1}{(1-z)^2} \prod^r_{i = 1} \frac{d_i z}{1 + (d_i - 1) z} .
\end{equation}
\label{lem:hirzebruch}
\end{lemma}
We give a proof for completeness:
\begin{proof}
Denote $X = V(g_1,\ldots,g_r)$ the complete intersection of generic forms.
We have a short exact sequence 
$$
0 \longrightarrow \mathcal{T}_X \longrightarrow \mathcal{T}_{\mathbb{CP}^n}|_X \longrightarrow \mathcal{N}_{X/\mathbb{CP}^n} \longrightarrow 0
$$
and that  
$$\mathcal{N}_{X/\mathbb{CP}^n} = \mathcal{O}_{\mathbb{CP}^n}(X)|_X = \bigoplus \mathcal{O}_X(d_i).$$
The Chern class of $X$ is thus
$$c(\mathcal{T}_X) = \frac{(1+h_X)^{n+1}}{\prod (1 + d_i h_X)},$$
where $h_X$ is the restriction of the hyperplane class in $\CP^n$ to $X$. 
The Euler characteristic of $X$ is the degree of $$c(\mathcal{T}_X) \cap [X] = \frac{(1+h)^{n+1}}{\prod (1 + d_i h)} \prod{d_i h} = (1+h)^{n+1} \prod \frac{d_i h}{1+d_i h}\in A^\ast (\CP^n),$$ see e.g. \cite[Chapter 3]{griffiths2014principles}. 
As $A^* \CP^n$ is a truncated polynomial ring in $h$, this degree is equal to the coefficient of $z^{n}$ in $(1+z)^{n+1} \prod \frac{d_i z}{1+d_i z} \in \C((z))$.
In other words, 
$$\chi(X) = \res_0 \frac{(1+z)^{n+1}}{z^{n+1}} \prod \frac{d_i z}{1+d_i z} dz.$$
Here, $\res_0 \sum_{i \in \Z} a_i z^i dz = a_{-1}$, the (formal) residue around 0. 

Consider the map $g(z) = \frac{z}{1-z}$. 
Note that  $(\frac{1+g(z)}{g(z)})^{n+1} = \frac{1}{z^{n+1}}$. 
Formal change of coordinates shows that $\res_0 f(z) d z = \res_0 \frac{f(g(z))}{(1-z)^2} d z$ for $f \in \C((z))$. 
In our setting, this gives
\begin{align*}
        \res_0 \frac{(1+z)^{n+1}}{z^{n+1}} \prod \frac{d_i z}{1+d_i z} dz &=
        \res_0 \frac{1}{(1-z)^2 z^{n+1}} \prod{\frac{d_i\frac{z}{1-z}}{1+{d_i}\frac{z}{1-z}}} dz  \\ &=
        \res_0 \frac{1}{(1-z)^2 z^{n+1}} \prod{\frac{d_i z}{1+(d_i-1)}} dz
\end{align*}
This is the same as the coefficient of $z^n$ in $\frac{1}{(1-z)^2} \prod{\frac{d_i z}{1+(d_i-1)}}$.
\end{proof}

For $r = 0$, the statement of the lemma reduces to $\chi(\CP^n) = n + 1$, which equals the coefficient of $z^n$ in $\frac{1}{(1-z)^2}$.

\begin{lemma}\label{lem: bn}
    Let $d_1, \dots d_r \geq 1$.
    Define $\tilde{b}_n[d_1, \dots d_r]$ via
    \begin{equation}
        \sum_{n \geq 0} \tilde{b}_n[d_1, \dots, d_r] z^n = \frac{1}{(1-z)^2} \prod^r_{i = 1} \frac{1- z}{1 + (d_i - 1) z} .
    \end{equation}
    Then for generic forms $g_1, \dots , g_r$ of degree $d_1, \dots d_r$ on $\CP^n$, we have that 
    $$\chi\big(\CP^n \setminus (V(g_1) \cup \dots \cup V(g_r))\big) = \tilde{b}_n[d_1, \dots d_r].$$
\end{lemma}
\begin{proof}
Let $[r] = \set{1, \dots, r}$. For $S = \set{i_1 < \dots < i_s} \subset [r]$,
let $c_n[S] = c_n[d_{i_1}, \dots, d_{i_s}]$. Then
we have that $\chi(V(\set{ g_i \colon i \in S})) = c_n[S]$ for $g_1, \dots, g_r$ in general position.

By inclusion-exclusion,
\begin{equation}
    \chi(\CP^n \setminus (V(g_1) \cup \dots \cup V(g_r))) = \sum_{S \subset [r]} (-1)^{\# S} V(\set{ g_i \colon i \in S}) = \sum_{S \subset [r]} (-1)^{\# S} c_n[S]
\end{equation}

Now
\begin{align*}
    \sum_{n \geq 0} \sum_{S \subset [r]} (-1)^{\# S} c_n[S] z^n &= 
    \sum_{S \subset [r]} (-1)^{\# S}  \frac{1}{(1-z)^2} \prod_{i \in S} \frac{d_i z}{1 + (d_i - 1) z} \\&=
    \frac{1}{(1-z)^2} \prod^r_{i = 1} \left( 1 - \frac{d_i z}{1 + (d_i - 1) z}\right) \\&=
    \frac{1}{(1-z)^2} \prod^r_{i = 1} \frac{1 -  z}{1 + (d_i - 1) z}
\end{align*}
This finishes the proof.
\end{proof}
\begin{remark}
We denote $\tilde{b}_n[1,d_1,\dots d_r] = b_n[d_1, \dots, d_r]$.
Then $b_i[d_1, \dots, d_r]$ is the coefficient of $z^i$ in  $$\frac{1}{(1-z)} \prod^r_{i = 1} \frac{1- z}{1 + (d_i - 1) z}.$$
And for generic forms $g_1, \dots, g_r$ on $\C^n$ of degree $d_1, \dots, d_r$, we have
$$
\chi(\C^n\setminus V(\prod_{i=1}^r g_i) = b_n[d_1,\ldots,d_r] .
$$
\end{remark}

\subsection{Proof of the main results}
In this section, we prove our main results for adding generic hypersurfaces to fixed hyperplane arrangements.

We start with the generic result Theorem \ref{thm: hyperplane multiple hypersurface}.
% \newtheorem*{restate}{Theorem \ref{thm: hyperplane multiple hypersurface}}
% \begin{restate}
%     Let $\A = \{V(f_1),\ldots,V(f_k)\} \subset \C^n$ be a hyperplane arrangement with characteristic polynomial $\chi_A(t) = \sum^n_{i=0} a_i t^i$. Let $d_1, \dots d_r \in \Z_+$.
%     Denote $b_i[d_1, \dots, d_r]$ to be the coefficient of $z^i$ in  $$\frac{1}{(1-z)} \prod^r_{i = 1} \frac{1- z}{1 + (d_i - 1) z}.$$
%     Then, for generic forms $g_1, \dots, g_r$ on $\C^n$ of degree $d_1, \dots, d_r$, we have: 
%     $$\chi(\C^n \setminus (\A  \cup V(\prod_{i=1}^rg_i))) = \sum^n_{i=0} a_i b_i[d_1, \dots, d_r].$$
% \end{restate}
\begin{proof}[Proof of Theorem \ref{thm: hyperplane multiple hypersurface}]
By inclusion-exclusion of the Euler characteristic, we have
\begin{equation*}
    \chi(\C^n \setminus (\A  \cup V(g_1 \cdots g_r))) = \sum_{\B \subset \A} (-1)^{\# \B} \chi( \bigcap \B \setminus V(g_1 \cdots g_r)).
\end{equation*}
Note that if $\bigcap \B$ is empty, then $\bigcap \B \setminus V(g_1 \cdots g_r)$ is empty as well and we have $\chi( \bigcap \B \setminus V(g_1 \cdots g_r))=0$. 
So, it is enough to sum over $\bigcap \B$ non-empty.

Consider the compactification $\C^n \subset \CP^n$ with $g_0$ the linear form vanishing on $\CP^n \setminus \C^n$. 
For $\B \subset \A$, we have that $\bigcap \B = \overline{\bigcap \B} \setminus V(g_0)$, so 
$$ \bigcap \B \setminus V(g_1 \cdots g_r) = \overline{\bigcap \B} \setminus V(g_0 g_1 \cdots g_r).$$
Let ${}^hg_1,\ldots,{}^hg_r$ be homogenization of $g_1,\ldots,g_r$ in $\CP^n$.
Since $g_1,\ldots,g_r$ are generic forms on $\C^n$, $g_0,{}^hg_1,\ldots,{}^hg_r$ are in general positions relative to $\overline{\bigcap \B} \cong \CP^{\dim \bigcap \B}$.
We can then apply Lemma \ref{lem: hirzeburch} to obtain 
$$
\chi(\overline{\bigcap \B} \setminus V(g_0 g_1 \cdots g_r)) = \tilde{b}_{\dim \bigcap \B}[1,d_1,\dots d_r] = b_{\dim \bigcap \B}[d_1, \dots, d_r].
$$

Together, this gives
\begin{equation}
    \chi(\C^n \setminus (\A  \cup V(g_1 \cdots g_r)))
    = \sum_{\substack{\B \subset \A\\ \bigcap \B \text{ non-empty}}} (-1)^{\# \B} b_{\dim \bigcap \B}[d_1, \dots, d_r]
    \label{eqn:result}
\end{equation}
The Whitney's Theorem presented in Proposition \ref{thm:whitney} above gives
\begin{equation}
\sum_{\substack{\B \subset \A\\ \B \text{ central}}} (-1)^{\# \B} z^{\dim \bigcap \B} = \chi_\A(z)
\label{eqn: whitney theorem}
\end{equation}

The result now follows from \cref{eqn:result} and \cref{eqn: whitney theorem}.
\end{proof}

When there is one generic hypersurface $V(g)$ of degree $d$ added, the above theorem simplifies as follows.

\begin{corollary}\label{cor: hyperplane one generic hypersurface}
Let $\A = \{V(f_1),\ldots,V(f_k)\} \subset \C^n$ be a hyperplane arrangement with characteristic polynomial $\chi_A(t)$ and $g$ is a generic degree $d$ form on $\C^n$. Then 
$$\chi\big(\C^n \setminus (\A  \cup V(g))\big) = \chi_\Acal(1-d).$$
\end{corollary}
\begin{proof}
% [Proof of Corollary \ref{cor: hyperplane one generic hypersurface}]
From Theorem \ref{thm: hyperplane multiple hypersurface}, 
$$
\chi\big(\C^n \setminus (\A\cup V(g))\big) = \sum_{i=0}^n a_i b_i[d],
$$
where $b_i[d]$ is the coefficient of $z^i$ in  $$\frac{1}{(1-z)} \frac{1- z}{1 + (d - 1) z} = \frac{1}{1+(d-1)z}.$$
Hence, 
$$
\chi\big(\C^n \setminus (\A\cup V(g))\big) = \sum_{i=0}^n a_i (1-d)^i = \chi_\A(1-d).
$$
\end{proof}

Theorem \ref{thm: one real point per region} is proved using Corollary \ref{cor: hyperplane one generic hypersurface} with a generic degree two hypersurface added.

% \newtheorem*{restateSECOND}{Theorem \ref{thm: one real point per region}}
% \begin{restateSECOND}
% Suppose $\Acal = \{V(f_1),\ldots,V(f_k)\}$ is an essential real hyperplane arrangement, i.e. there exists some $S\subseteq [k]$ with $\bigcap_{i\in S} V(f_i)$ is a point.
% Let $g$ be a generic degree 2 polynomial satisfying $g(x)>0$ for all $x\in \R^n$. Then for the function
% $$
% \psi(x) = \sum_{i=1}^k u_i \log |f_i(x)| -v \log|g(x)|
% $$
% where $u_1,\ldots,u_k,v\in \R_{>0}$ and $2v>\sum_{i=1}^k u_i$, all of its critical points are real and there is precisely one critical point per region for $\R^n \setminus \Acal$.
% \end{restateSECOND}
\begin{proof}[Proof of Theorem \ref{thm: one real point per region}]
The hyperplane arrangement $\Acal = \{V(f_1),\ldots,V(f_k)\}$ is essential, so $\C^n\setminus \Acal$ is very affine. It follows that $\C^n\setminus (\Acal \cup V(g) )$ is also very affine and by \cite{huh2013maximum}, the number of critical points of $\psi(x_1,\ldots,x_n)$ equals $\chi\big( \C^n\setminus (\Acal \cup V(g) ) \big).$
By Corollary~\ref{cor: hyperplane one generic hypersurface}, $\chi\big( \C^n\setminus (\Acal \cup V(g) ) \big) = \chi_\Acal(-1)$, which coincides with the number of regions for $\R^n\setminus \Acal$. 

When we choose $u_1,\ldots,u_k,v\in \R_{>0}$ and $2v>\sum_{i=1}^k u_i$, $\psi(x)$ tends to negative infinity at the boundary of each region and infinity in $\R^n\setminus\Acal$ ensuring at least one local maximum of $\psi(x)$ in each region.
The total number of complex critical points of $\psi(x)$ equals the number of regions. 
So, we must have all complex critical points of $\psi$ are real and there is exactly one real critical point per region of $\R^n\setminus\Acal$.
\end{proof}

\subsection{Algorithm for sampling one interior point per region in hyperplane arrangement complement}
We present in this subsection an algorithm based on Theorem \ref{thm: one real point per region} to sample one interior point per region of $\R^n\setminus \Acal$, where $\Acal = \{V(f_1),\ldots, V(f_k)\}$ is a real hyperplane arrangement.
We run numerical experiments and compare the performance of our algorithm with other existing methods in terms of run-time. The code for the experiments can be found at \url{https://github.com/QWE123665/hypersurfacearrangements}.

The main task of our algorithm is to solve critical points of $\psi(x) = \sum_{i=1}^k u_i \log |f_i(x)| -v \log|g(x)|$, where $g$ is a generic degree two polynomial satisfying $g(x)>0$ for all $x\in \R^n$ and $u_1,\ldots,u_k,v\in \R_{>0}$ and $2v>\sum_{i=1}^k u_i$.
Leveraging advances in numerical algebraic geometry, this task is efficiently performed using monodromy and homotopy continuation methods, as implemented in the Julia package \texttt{HomotopyContinuation.jl} \cite{breiding2018homotopycontinuation}.

First, we construct a start system by treating $u_1,\ldots,u_k,v$ as unknown parameters and run monodromy solve.
The stopping criterion is that the total number of solutions equals the number of regions in $\R^n\setminus\Acal$, which can be computed via the Julia package \texttt{CountingChambers.jl} in \cite{brysiewicz2023computing} if the coefficients of the hyperplanes are rational. 
The start system yields all the solutions of $\psi(x)$ for one specific parameter choice $(u_1,\ldots,u_k,v)$.
Next, we track these solutions to the system corresponding to our desired parameter values satisfying $u_1,\ldots,u_k,v\in \R_{>0}$ and $2v>\sum_{i=1}^k u_i$. 
Each solution for this system lies in a distinct region of $\R^n\setminus\Acal$, verified by checking the sign patterns of $(f_1,\ldots,f_k)$ at the solutions are distinct. 
The algorithm is summarized in Algorithm \ref{alg:morse}.

\begin{algorithm}[ht]
%\scriptsize
\caption{Sampling points in each region of $\R^n\setminus \Acal$ (Morse)} \label{alg:morse}
\KwIn{$f_1,\ldots,f_k \in \Q[x_1,\ldots,x_n]$ where $\Acal = \{V(f_1),\ldots, V(f_k)\}$}

\KwOut{A list of points $\{p_1,\ldots,p_N\}$ with one point in each region of $\R^n \setminus \Acal$}

Choose $g(x) = \sum_{i=1}^n (x_i-a_i)^2 + (\sum_{i=1}^n b_ix_i)^2 +1 $ for random choices of $a_1,\ldots,a_n, b_1,\ldots,b_n$.

Compute the number of regions $N$ of $\R^n\setminus \Acal$ using \texttt{CountingChambers.jl}.

Use monodromy solve implemented in \texttt{HomotopyContinuation.jl} to compute the critical points of $\psi(x) = \sum_{i=1}^k u_i \log |f_i(x)| -v \log|g(x)|$ with parameters $u_1,\ldots,u_k,v$. 
The stopping criterion is \verb|target_solutions_count = N|.

Track the system obtained to the system with parameters $(u_1,\ldots,u_k,v) = (1,\ldots,1,\lceil \frac{k+1}{2}\rceil )$ and obtain solutions $p_1,\ldots,p_N$.

Evaluate $(f_1,\ldots,f_k)$ at $p_1,\ldots,p_N$ and check the signs of the evaluations are distinct.
\end{algorithm}

We compare our algorithm to three other methods.

Denote $H_+^i = \{x: f_i(x)>0\}$ and $H_\_^i = \{x: f_i(x)<0\}$ for $i=1,\ldots,k$.
The first two methods run through all possible intersections of the half spaces $V = V_1\cap \ldots\cap V_k$ where $V_i\in \{H_+^i,H_\_^i\}$, which is either empty or a polyhedron.

For each intersection of the half-spaces $V = V_1\cap \ldots \cap V_k$, the first method computes its dimension using the software \texttt{Polymake} \cite{gawrilow2000polymake}. If the dimension of $V$ is $n$, we compute the V-representation of the polyhedron from its H-representation. 
The information of the V-representation then gives an interior point of the region.
We denote this method (BFHV).

The second method uses linear programming.
For each intersection of the half-spaces $V = V_1 \cap \ldots \cap V_k$, we define a linear programming problem with constraints taken from $V_1,\ldots, V_k$ and objective function zero.
Via the interior-point-method \cite{karmarkar1984new}, solving the linear programming returns a point in $V$.
We denote this method (BFIPM).

The third method is from \cite{KastnerPanizzutHyperplanes}, where the authors compute cell decompositions for central hyperplane arrangements with specified supporting cones by computing full-dimensional cones of certain normal fans.
Unlike the first two methods, this approach avoids exhaustive enumeration and instead flips hyperplanes to determine neighboring cones. 
Any affine hyperplane arrangement can be handled by computing the cell decomposition of the central hyperplane arrangement of its homogenization and using the supporting cone $\{x_0>0\}$.
Interior points are obtained as positive linear combinations of the rays of each cell. 
We refer to this method as (CD).

\begin{figure}[htbp]
    \centering
    \includegraphics[scale = 0.8]{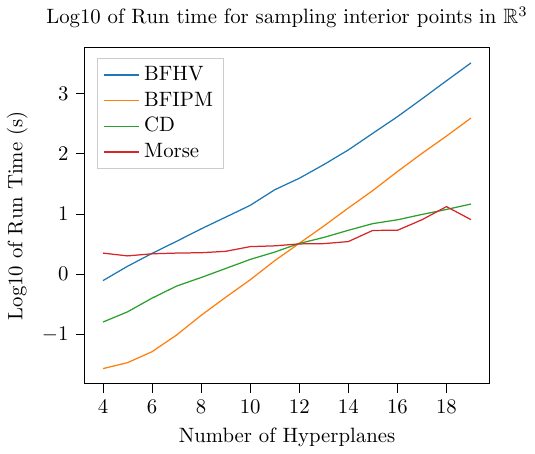}
    \includegraphics[scale = 0.8]{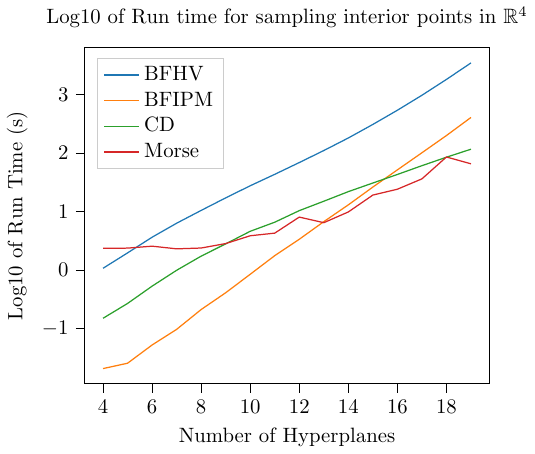}
    \caption{The time in log-scale of sampling one interior point per region for 4 to 19 random hyperplanes in $\R^3$(left) and $\R^4$(right) for the four methods.}
    \label{fig:R^3R^4randomhyperplanes}
\end{figure}

We compare the run times of our method to the other three methods for sampling one interior point per region for 4 to 19 hyperplanes in $\R^3$ and $\R^4$.
The hyperplanes have random coefficients chosen as integers between -10 and 10 across the four methods.
The results are shown in Figure \ref{fig:R^3R^4randomhyperplanes}. 
It can be seen that our algorithm Morse is comparable to CD and significantly faster than BFIPM and BFHV.

\begin{figure}[htbp]
    \centering
    \includegraphics[scale = 0.8]{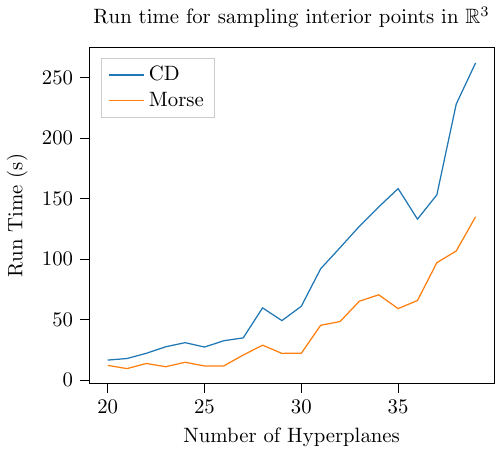}
    \includegraphics[scale = 0.8]{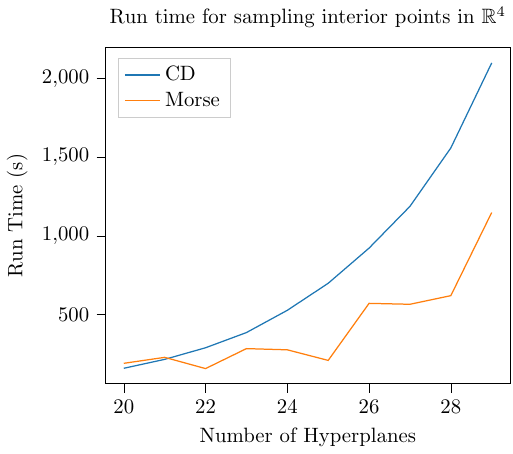}
    \caption{The time in log-scale of sampling one interior point per region for 20 to 39 random hyperplanes in $\R^3$(left) and 20 to 29 in $\R^4$(right) for CD and Morse.}
    \label{fig:morehyperplanes}
\end{figure}

We compare the run times of our algorithm with CD for the task of sampling one interior point per region for 20 to 39 hyperplanes in $\R^3$ and 20 to 29 hyperplanes in $\R^4$.
The hyperplanes have random coefficients chosen as integers between 
-10 and 10.
The results are shown in Figure \ref{fig:morehyperplanes}.
The Morse method is significantly faster than the CD method, with runtime improvements that become more pronounced as the number of hyperplanes increases.

We apply the algorithms to the resonance arrangements discussed in the introduction.
In other words, for each dimension $d$, we take $2^d-1$ hyperplanes normal to vectors with entries 0 or 1 but not all 0. We compare the run-time of the four methods applied to $d=2,\ldots,6$. The results are shown in Table \ref{table:hyperplanes normal to 0 1}. It can be seen that our method Morse is comparable to the other methods for $d=2,3$ and significantly faster for $d=4,5,6$.

\begin{table}[htbp] 
\centering
\begin{tabular}{|c|c|c|c|c|c|c|} 
  \hline
  d & $\#$ hyperplanes & $\#$ regions & BFHV(s) & BFIPM(s) & CD(s) &  Morse(s)
 \\
 \hline
 2 & 3 & 6 & 0.27 & 0.02 & 0.10 & 1.33
 \\
  \hline 
  3& 7 & 32 & 4.38 & 0.12 & 0.40 &1.40
   \\
  \hline 
  4 & 15 & 370 & 1235.4 &32.5 & 4.98&1.84
   \\
  \hline 
  5 & 31 & 11292 & - & -& 209.19& 34.44
   \\
  \hline 
  6 & 63 & 1066044& -& -& 40517.84&4557.80
  \\
  \hline
\end{tabular}
\caption{Run-time for sampling interior points from the complement of hyperplanes normal to vectors with entries 0 or 1. The empty entries mean that the experiment does not terminate in 24h.}
\label{table:hyperplanes normal to 0 1}
\end{table}

\section{Hypersurface arrangement with generic hypersurfaces added}\label{sec: hypersurface}
\subsection{Preliminaries}\label{sec: csm and milnor numbers}
In this subsection, we briefly introduce the Milnor numbers of hypersurfaces and the Chern-Schwartz-MacPherson (CSM) class, which prepares us for the proof of Theorem \ref{thm: hypersurface arrangement with multiple hypersurface} in the next subsection.

Consider a holomorphic function germ $f: (\C^n,0) \to (\C,0)$, the \emph{Milnor number} of $f$, introduced by John Milnor \cite{milnor1968singular}, is defined to be $\dim_\C \C\{x_1,\ldots,x_n\}/Jf$, where $Jf$ is the ideal generated by partial derivatives of $f$. It coincides with the topological degree of the map $\set{z \in \C^n | \| z \| = \epsilon} \rightarrow \set{z \in \C^n | \| z \| = 1}, z\mapsto \frac{\nabla f}{\|\nabla f\|}$, where $\epsilon>0$ is small and $\nabla f$ is the gradient of $f$. 

When $f$ is a homogeneous polynomial defined on $\CP^n$, we denote $\mu^i(f)$ to be the topological degree of $\nabla f|_{H_i}$, the polar map of the restriction of $f$ to a generic linear subspace $H_i$ of dimension $i$.
%$$
%\nabla f: \CP^n \dashrightarrow \CP^n, \text{where } \nabla f(x) = %[\frac{\partial f}{\partial x_0}(x),\ldots,\frac{\partial f}{\partial x_n}%(x)]
%$$
%restricted to a generic $\CP^i\subseteq \CP^n$. \br{I think it should be }
Teissier \cite{teissier1973cycles} relates the Milnor numbers $\{\mu^i(f)\}$ with the mixed-multiplicities of the $i$-th mixed multiplicity of the ideals $(x_0,\ldots,x_n)$ and $(\frac{\partial f}{\partial x_0},\ldots,\frac{\partial f}{\partial x_n})$.
Another interpretation is that the graph $\Gamma_f$ of the polar map $\nabla f$ has class 
$$
[\Gamma_f]= \sum_{i=1}^n \mu^i(f)[\CP^{n-i}\times \CP^{i}]
$$
in the Chow group $A_n(\CP^n\times \CP^n)$.

The Milnor numbers $\{\mu^i(f)\}_{i=0}^n$ also have the following topological interpretation:
denote $D(f)_i$ (resp. $V(f)_i$ ) to be the space $\{x\in \CP^n: f(x)\neq 0\}$ (reps. $\{x\in \CP^n: f(x)= 0\}$ )
restricted to a generic $\CP^i\subseteq \CP^n$.
Dimca and Papadima \cite{dimca2003hypersurface} showed that 
$D(h)_i$ is homotopy equivalent to $D(h)_{i-1}$ with $\mu^i(f)$ dimension-$i$ cells attached. 
It immediately provides the formula for calculating the Euler characteristic of $D(f)$ using Milnor numbers.

\begin{lemma}\cite{dimca2003hypersurface}
Suppose $f$ is a homogeneous polynomial defining a hypersurface in $\CP^n$. Then,
\begin{equation}\label{eq: euler via milnor numbers}
\chi_\C(D(f)) = \sum_{i=1}^n (-1)^i \mu^i(f).
\end{equation}
\end{lemma}

The characteristic polynomial of a hyperplane arrangement relates to the Milnor numbers of the product of the corresponding linear forms in the following way.

\begin{proposition}\label{thm: milnor numbers and characteristic poly of hyperplane}
Suppose $f$ is a product of homogeneous linear forms. Let $\widetilde{\Acal}$ be the central hyperplane arrangement in $\C^{n+1}$ defined by linear forms in $f$ and let $\Acal^H$ be the hyperplane arrangement in $\C^n$ by declaring $H\in \widetilde{\Acal}$ to be the hyperplane at infinity. Then,
\begin{enumerate}[(i)]
\item \cite{StanleyHyperplanes} We obtain $\chi_{\widetilde{\Acal}}(z)/(z-1) = \chi_{\Acal^H}(z)$.
\item \cite[Corollary 25]{Huh2012} We have
$$\chi_{\Acal^H}(z) = \sum_{i=0}^n (-1)^i \mu^i(f) z^{n-i}$$ 
and $\{\mu^i(f)\}_{i=0}^n$ coincides with the Betti-numbers of $D(f)$.
\end{enumerate}
\end{proposition}

Now we turn to the CSM class \cite{aluffi2005characteristic}. The CSM class is a generalization of the Chern class to singular varieties. 
For a variety $X$, one can define the CSM class of any constructible function of the form $\sum_{i=1}^\ell n_i \mathds{1}_{Z_i}$ where $n_i\in \Z$ each $Z_i$ is a subvariety of $X$. 
When $X$ is nonsingular and complete, the CSM class of $X$ denoted by $c_{SM}(\mathds{1}_X)$ or simply 
$c_{SM}(X)$ coincides with $c(\mathcal{T}_X)\cap [X]$ where $c(\mathcal{T}_X)$ is the Chern class of the tangent bundle and $[X]$ is the fundamental class for $X$ in the appropriate Chow group. We list three fundamental properties of the CSM class:
\begin{enumerate}
    \item The degree of $c_{SM}(X)$ agrees with $\chi(X)$.
    \item $f_\ast(c_{SM}(\psi)) = c_{SM}(f_\ast(\psi))$ for $f: X\to Y$ proper and $\psi: X\to \Z$ constructible.
    \item If $X_1,X_2\subseteq X$, then $\csm(X_1\cup X_2 ) = \csm(X_1) + \csm(X_2)-\csm(X_1\cap X_2)$ in $A_\ast(X)$. 
\end{enumerate}

The Milnor numbers of a hypersurface in $\CP^n$ relate to the CSM class in the following way.
\begin{theorem}\cite[Theorem 2.1]{aluffi2003computing}\label{thm: csm and milnor}
Suppose $f$ is a homogeneous polynomial defining a hypersurface in $\CP^n$ and denote $D(f)$ as the hypersurface's complement.
Then,
$$
\csm(D(f)) = \sum_{i=0}^n (-1)^i \mu^i(f) h^i(1+h)^{n-i} \in A_\ast(\CP^n),
$$
where $\{\mu^i(f)\}_{i=0}^n$ are Milnor numbers and $h$ is the class of a hyperplane in $A_\ast(\CP^n)$.
\end{theorem}

We conclude this subsection with an example.
\begin{example}
Consider the hyperplane arrangement $\Acal$ of three coordinate planes defined by $f = xyz$ in $\CP^2$. We first compute the Milnor numbers of $f$.

The polar map $[x:y:z]\mapsto [yz:xz:xy]$ is the standard Cremona birational involution, so in particular it has topological degree 1. So, $\mu^2(f) = 1$.
For $\mu^1(f)$, we consider the polar map of the restriction to $H:ax+by+cz =0$ with $a,b,c\in \C$ generic. Without loss of generality, suppose $c=1$.
Then $f|_H = xy(-ax-by)$, so the polar map of the restriction is equivalent to $[x:y]\mapsto [2axy+by^2:ax^2+2bxy].$  
Solving $2axy+by^2 = \alpha, ax^2+2bxy = \beta $ for generic $\alpha,\beta \in \C$ gives four solutions in $\C^2$ and two solutions in $\CP^1$, hence $\mu^1(f)=2$.
The Milnor number $\mu^0(f)$ is the topological degree of the polar map restricted to a point, so $\mu^0(f) = 1$.

We can choose $z=0$ as the hyperplane at infinity and obtain a decone $\Acal^H \subseteq \C^2$ of $\Acal$ consisting of two lines $x = y =0$. The characteristic polynomial $\chi_{\Acal^H}(z) = z^2-2z+1$ which coincides with $\mu^0(f)z^2 - \mu^1(f)z +\mu^2(f)$, confirming Proposition \ref{thm: milnor numbers and characteristic poly of hyperplane}(b).

We also compute the CSM class of $D(f)$ by considering the Euler characteristic of its intersection with generic hyperplanes.
The variety $D(f)$ is homotopy equivalent to $\C^2$ taking out two lines that intersect at a point. By inclusion-exclusion formula, $\chi(D(f)) = 0$. Let $H$ be a generic hyperplane in $\CP^2$, then $D(f)\cap H$ is homotopy equivalent to a line in $\C^2$ taking out two points. So, $\chi(D(f)\cap H) = -1$.
Let $H' \subseteq \CP^2$ be another generic hyperplane, then $D(f)\cap H \cap H'$ is a single point so $\chi(D(f)\cap H\cap H') = 1$. 
Denote $\chi_{D(f)}(z) = \chi(D(f)\cap H\cap H') (-z)^2 + \chi(D(f)\cap H)(-z) + \chi(D(f)) (-z)^0 = z^2 + z.$
By \cite[Theorem 1.1]{Aluffi2013}, the CSM class of $D(f)$ is 
$\frac{z \cdot \chi_{D(f)}(-1-z)+ \chi_{D(f)}(0)}{z+1} = z^2$ in $A_\ast(\CP^2)$ where $z^i$ represents the class $[\CP^i]$.

Using Theorem \ref{thm: csm and milnor}, we have $\csm(D(f)) = h^2 - 2h(1+h) + (1+h)^2 = 1 $, where $H$ represents the class of a hyperplane. 
In other words, we obtain $\csm(D(f))= [\CP^2]$, obtaining the same answer as above.
\end{example}

\subsection{Prove of Theorem \ref{thm: hypersurface arrangement with multiple hypersurface}}

In this subsection, we prove Theorem \ref{thm: hypersurface arrangement with multiple hypersurface}, a generalization of Theorem \ref{thm: hyperplane multiple hypersurface} for the hyperplane arrangement to the hypersurface arrangement. 
As seen in Section \ref{sec: csm and milnor numbers}, the coefficients of the characteristic polynomial of hyperplane arrangements generalize to Milnor numbers for hypersurface arrangements, and the Milnor numbers closely relate to the CSM class of the hypersurface. 
Hence, we start by understanding how the CSM class of some variety $X$ changes when its intersection with a generic hypersurface is removed.

\begin{theorem}\label{thm: csm class after intersection}
Suppose $V$ is a non-singular projective variety and $X\subseteq V$ is a hypersurface. Suppose $G$ is a generic hypersurface of degree $d$ in $V$, then 
\begin{equation}\label{eq: csm class add hypersurface}
c_{SM}(X\setminus G) = \frac{1}{1+G} \cdot c_{SM}(X)
\end{equation}
in $A_\ast(X).$
\end{theorem}
\begin{proof}
Equation \eqref{eq: csm class add hypersurface} is equivalent to 
$c_{SM}(X\cap G) = \frac{G}{1+G} \cdot c_{SM}(X)$.
By the proof of \cite[Theorem 3.1]{aluffi2011chern}, we need to show that $$s(\Sing (X\cap G), G) = G \cdot s(\Sing (X), V)$$ in $A_\ast(X)$, where $s(\cdot,\cdot)$ denotes the Segre class and $\Sing(\cdot)$ denotes the singular locus of a variety.

Let $\phi: \tilde{V}\rightarrow V$ be the blow up of $V$ along $\Sing X$.
By the proof of \cite[Claim 3.2]{aluffi2011chern}, if $\phi^{-1}(G)$ is the same as the strict transform $\tilde{G}$ of $G$, then 
$$s(\Sing (X)\cap G, G) = G \cdot s(\Sing (X), V).$$
Let $\Ical,\Jcal$ be the ideal sheaves of $\Sing X$ and $G$ in $V$.
By definition, $\Proj(\bigoplus_j \Ical^j) = \tilde{V}$.
We also have $\phi^{-1}(G) = \Proj(\bigoplus_j \frac{\Ical^j}{\Jcal \Ical^j})$ and 
$\tilde{G} = \Proj(\bigoplus_j \frac{(\Ical+\Jcal)^j}{\Jcal}) = \Proj(\bigoplus_j \frac{\Ical^j+\Jcal}{\Jcal})$. 
To show $\phi^{-1}(G) = \tilde{G}$, it is enough to show that $\Jcal \cap \Ical^j = \Jcal \Ical^j$ for $j>>0$ since then 
$\frac{\Ical^j}{\Jcal \Ical^j} = \frac{\Ical^j+\Jcal}{\Jcal}$ for $j>>0$ and hence $\Proj(\bigoplus_j \frac{\Ical^j}{\Jcal \Ical^j})=\Proj(\bigoplus_j \frac{\Ical^j+\Jcal}{\Jcal})$. 
$V$ is a projective variety, so it can be covered by finitely many affine varieties. 
To show $\Jcal \cap \Ical^j = \Jcal \Ical^j$ for $j>>0$, it suffices to show that in the coordinate ring of each affine variety, the ideals $I, J$ corresponding to the ideal sheaves $\Ical, \Jcal$, for $j>>0$ satisfy $I^j\cap J = I^j J$. 
Indeed, this is true for any $j>0$.
It is obvious that $I^j J \subseteq I^j \cap J$.
Since $G$ is a generic hypersurface, $J=\langle g \rangle $ for some $g$ irreducible and $g\notin I^j$. Then $J+I^j = \langle 1 \rangle$. 
Suppose $a \in J, b\in I^j$ with $a+b = 1$ and $c \in I^j \cap J$. 
Then, $c = c(a+b) = ca+cb$ and $ca,cb \in I^j J$ since $c \in I^j \cap J, a\in J, b\in I^j$. Hence, $c\in I^j J$ and $I^j J \supseteq I^j\cap J$.

We are left to show that 
$$s(\Sing (X\cap G), G) = s(\Sing X \cap G, G)$$ in $A_\ast(X)$. By Teissier’s idealistic Bertini \cite[Section 2.8]{teissier1977hunting}, for generic hypersurface $G$, the ideals defining $\Sing (X\cap G)$ and $\Sing (X) \cap G$ have the same integral closure. 
Segre classes are birational invariant, hence $s(\Sing (X\cap G), G) = s(\Sing (X) \cap G, G)$ by considering their normalized blow-ups.
\end{proof}

\begin{corollary}\label{cor: csm remove hypersurface intersection}
Let $U\subseteq \CP^n$ be the complement of the union of hypersurfaces $
\bigcup_{i=1}^k V(f_i)$. 
Let $G$ be a generic hypersurface of degree $d$, then 
$$
c_{SM}(U \setminus G)= \frac{1}{1+dh}\cdot c_{SM}(U)
$$
in $A_\ast(\CP^n)$, where $h$ is the class of a hyperplane.
\end{corollary}
\begin{proof}
By the inclusion-exclusion formula for the CSM classes, it suffices to show that for any $S \subseteq [k]$, 
$$
c_{SM}(V(\prod_{i\in S} f_i) \setminus G)= \frac{1}{1+dh}c_{SM}(V(\prod_{i\in S} f_i) )
$$
in $A_\ast(\CP^n).$

Denote $X = V(\prod_{i\in S} f_i)$. 
By Theorem \ref{thm: csm class after intersection},
$c_{SM}(X \cap G)= \frac{1}{1+G}\cdot c_{SM}(X)$ in $A_\ast(X)$.
Denote $i: X\hookrightarrow \CP^n$ the inclusion map. 
Then applying the pushforward map $i_\ast$, we have $i_\ast(c_{SM}(X \cap G))=c_{SM}(X \cap G),i_\ast(c_{SM}(X))=c_{SM}(X)$ and $i_\ast(G)=dh$. 
By the projection formula \cite[Theorem 3.2(c)]{fulton2013intersection}, we obtain 
$i_\ast(G^\ell \cdot c_{SM} (X)) = (dh)^\ell \cdot c_{SM}(X) $ for any $\ell\geq 0$.
It follows that $i_\ast(\frac{1}{1+G}\cdot c_{SM}(X))= \frac{1}{1+dh}\cdot c_{SM}(X)$ and hence $c_{SM}(X \setminus G)= \frac{1}{1+dh}\cdot c_{SM}(X)$ in $A_\ast(\CP^n).$
\end{proof}

\begin{corollary}\label{cor: csm remove multiple hypersurfaces}
Let $g_1,\ldots,g_r$ be generic forms on $\C^n$ with degrees $d_1,\ldots,d_r$.
Then for hypersurfaces $V(f_1),\ldots, V(f_k)$ with $\Acal = \bigcup_{i=1}^k V(f_i)$,
we have
$$
c_{SM}(\C^n\setminus (\Acal \cup \bigcup_{i=1}^r V(g_i) )) = \prod_{i=1}^r \frac{1}{1+d_i h} c_{SM}(\C^n \setminus \Acal).
$$
\end{corollary}
\begin{proof}
The result immediately follows from applying Corollary \ref{cor: csm remove hypersurface intersection} to $c_{SM}(\C^n\setminus (\Acal \cup V(g_1)),\ldots, c_{SM}(\C^n\setminus (\Acal\cup \bigcup_{i=1}^r V(g_i) ))$.
\end{proof}

We now prove Theorem \ref{thm: hypersurface arrangement with multiple hypersurface}. 
\begin{proof}[Proof of Theorem \ref{thm: hypersurface arrangement with multiple hypersurface}]
We have that 
\begin{align*}
c_{SM}(\C^n \setminus (\A  \cup V(g_1 \cdots g_r))) = & 
\left(\prod_{i=1}^r \frac{1}{1+d_i h} \right)c_{SM} (\C^n \setminus \Acal) & (\text{By Corollary \ref{cor: csm remove multiple hypersurfaces}})\\
=& \left(\prod_{i=1}^r \frac{1}{1+d_i h} \right) \left(\sum_{i=0}^n (-1)^i \mu^i({}^hfx_0) h^i (1+h)^{n-i} \right) & (\text{By Theorem \ref{thm: csm and milnor}}),
\end{align*}
where ${}^hf$ is the homogenization of $f$ using $x_0$ and $\mu^i(\cdot)$ denotes the $i$-th Milnor number of a hypersurface in $\CP^n$.
Therefore, $\chi(\C^n \setminus (\A  \cup V(g_1 \cdots g_r)))$ is the degree of $c_{SM}(\C^n \setminus (\A  \cup V(g_1 \cdots g_r)))$, and hence equals 
$$\text{the coefficient of $z^n$ in } \left(\prod_{i=1}^r\frac{1}{1+d_i z} \right) \left(\sum_{i=0}^n (-1)^i \mu^i(fx_0) z^i (1+z)^{n-i} \right).$$ 

Now, we are left to show that 
the coefficient of $z^{n-i}$ in $\left(\prod_{i=1}^r \frac{1}{1+d_i z}\right) (1+z)^{n-i}$ equals $b_{n-i}[d_1,\ldots,d_r]$.
Equivalently, we need to show
$$
\text{coeff of $z^{n-i}$ in } \left(\prod_{i=1}^r \frac{1}{1+d_i z}\right) (1+z)^{n-i} = \text{coeff of  $z^{n-i}$ in }
\frac{1}{(1-z)} \prod^r_{i = 1} \frac{1- z}{1 + (d_i - 1) z}.
$$
The coefficient of $z^{n-i}$ in $\left(\prod_{i=1}^r \frac{1}{1+d_i z}\right) (1+z)^{n-i}$
equals $\res_0 \left(\prod_{i=1}^r \frac{1}{1+d_i z}\right) \frac{(1+z)^{n-i}}{z^{n-i+1}} dz.$
Denote $q(z) =\left(\prod_{i=1}^r \frac{1}{1+d_i z}\right) \frac{(1+z)^{n-i}}{z^{n-i+1}}$ and $g(z) = \frac{z}{1-z}$.
Formal change of coordinates shows that $\res_0 q(z) dz = \res_0 \frac{q(g(z))}{(1-z)^2}$.
Note that
$$
\frac{q(g(z))}{(1-z)^2} =\left(\prod_{i=1}^r \frac{1-z}{1+(d_i-1)z}\right)\frac{1}{z^{n-i}} \frac{(1-z)}{z} \frac{1}{(1-z)^2}
= \frac{1}{(1-z)}\left(\prod_{i=1}^r \frac{1-z}{1+(d_i-1)z}\right)\frac{1}{z^{n-i+1}}.
$$
Hence, $ \res_0 \frac{q(g(z))}{(1-z)^2}$ agrees with the coefficient of $z^{n-i}$ in $\frac{1}{(1-z)}\left(\prod_{i=1}^r \frac{1-z}{1+(d_i-1)z}\right)$.
\end{proof}
\begin{remark}
Theorem \ref{thm: hyperplane multiple hypersurface} can alternatively be derived from the preceding theorem and Proposition \ref{thm: milnor numbers and characteristic poly of hyperplane}(ii). Nevertheless, we retain the current proof for its more combinatorial nature.
\end{remark}

\subsection{Upper bounds on the number of regions for hypersurface arrangements}

In this subsection, we consider real coefficients polynomials $f_1,\ldots,f_k$ and we study the number of regions (or connected components) for $\R^n\setminus \Acal$ where $\Acal = \cup_{i=1}^k V(f_i)$.

Suppose $g(x_1,\ldots,x_n)\in \R[x_1,\ldots,x_n]$ is a generic degree two polynomial such that $g(x)>0$ for any $x\in \R^n$. For generic $(u_1,\ldots,u_k,v)\in \R_{+}$ with $2v>\sum_{i=1}^k u_i \deg(f_i)$, the function 
$$
-\psi(x) = \sum_{i=1}^k u_i \log |f_i(x)| -v \log|g(x)|
$$ 
is Morse and it
goes to infinity at points on $\Acal$ and at infinity. So, each region of $\R^n\setminus \Acal$ contains at least one local minimum of $\psi$. This observation results in the proof of Theorem \ref{thm: number of regions bounded by complex euler}.

\begin{proof}[Proof of Theorem \ref{thm: number of regions bounded by complex euler}]
We use the notations defined in the paragraph above.
Since each region of $\R^n-\Acal$ contains at least one local minimum of $-\psi$, we have the number of regions for $\R^n \setminus \Acal$ is bounded above by the number of critical points of $\psi$. 
The assumption of $\R^n \setminus \Acal$ being smooth and very affine implies that $\R^n \setminus (\Acal\cup V(g))$ is also smooth and very affine. By \cite{huh2013maximum}, the number of  critical points of $-\psi$ coincides with $(-1)^n \chi(\C^n \setminus (\Acal\cup V(g))$. Applying Theorem \ref{thm: hypersurface arrangement with multiple hypersurface}, we have
$$
(-1)^n \chi\big(\C^n \setminus (\Acal\cup V(g))\big) = (-1)^n \sum_{i=0}^n (-1)^{n-i} \mu^{n-i}({}^hf x_0) b_i[2],
$$
where $b_i[2]$ is the coefficient of $z^i$ in $\frac{1}{1+z}$, so $b_i[2] = (-1)^i$.
We obtain 
$$(-1)^n\chi\big(\C^n \setminus (\Acal\cup V(g))\big) = \sum_{i=0}^n \mu^{n-i}({}^hf x_0)$$ and hence the first inequality in the theorem.

We are left to prove $\mu^{i}({}^hf x_0)\leq d^i$ for all $i=0,\ldots,n$. Recall that $\mu^{i}({}^hf x_0)$ is the topological degree of the polar map $[x_0,\ldots,x_n] \mapsto [\frac{\partial {}^hf x_0}{\partial x_0},\ldots, \frac{\partial {}^hf x_0}{\partial x_n}]$ restricted to a generic $\CP^i$, where each $\frac{\partial {}^hfx_0}{\partial x_i}$ has degree $d$. 
So, $\mu^{i}({}^h f x_0)$ is the number of intersection points for $i$ degree $d$ polynomials and by B\'ezout's Theorem, $\mu^{i}({}^h f x_0)\leq d^i$.
\end{proof}

\begin{remark}
The equality of $\# \text{ regions of } \R^n \setminus \Acal \leq \sum_{i=0}^n\mu^i(f^h)$ holds when all complex critical points of $-\psi(x)$ are real and there is exactly one critical point per region of $\R^n\setminus \Acal$.
One class of examples is when $\Acal$ is an essential hyperplane arrangement.
Another class arises from the study of Determinantal Point Processes from an algebraic statistics point of view \cite[Theorem 1.3]{friedman2024likelihood}, where $\Acal$ is a certain collection of hyperplanes and hypersurfaces of degree two.
\end{remark}

Using Morse theory, we can obtain a tighter inequality regarding the number of regions.
\begin{theorem}\label{thm: upper bound for number of regions using real euler char}
Suppose $f_1,\ldots,f_k\in \R[x_1,\ldots,x_n]$ and $\C^n\setminus(\Acal \cup V(g))$ is smooth and very affine where $\Acal = \cup_{i=1}^k V(f_i)$ and $g(x_1,\ldots,x_n)\in \R[x_1,\ldots,x_n]$ is a generic degree two polynomial such that $g(x)>0$ for any $x\in \R^n$. Then 
\begin{align*}
\# \text{ regions of } \R^n \setminus \Acal & \leq  \frac{1}{2}\bigl( |\chi\big(\C^n \setminus (\Acal\cup V(g))\big)| + \chi(\R^n \setminus \Acal) \bigr)\\
&=\frac{1}{2}\bigl(\sum_{i=0}^n\mu^i(f^hx_0)  + \chi(\R^n \setminus \Acal) \bigr).
\end{align*}
\end{theorem}
\begin{proof}
We continue to use the notations defined in the paragraph above Theorem \ref{thm: number of regions bounded by complex euler}. By the proof of Theorem \ref{thm: number of regions bounded by complex euler}, the total number of critical points of $-\psi(x)$ is $\sum_{i=0}^n\mu^i(f^hx_0)$. 
By genericity of $g(x)$ and $u_1,\ldots,u_k,v$ that defines $-\psi(x)$, $-\psi(x)$ is Morse and exhaustive on $\R^n\setminus \Acal$. In other words, all Hessian matrices of the real critical points for $-\psi(x)$ are nonsingular and all sublevel sets of $-\psi(x)$ are compact. 
In Morse theory (see e.g. \cite{nicolaescu2007invitation}), the index of a critical point for $-\psi(x)$ is defined to be the number of negative eigenvalues for the corresponding Hessian matrix and $\chi(\R^n\setminus \Acal)  = \sum_{i=0}^n (-1)^i \mu_i$ where $\mu_i$ is the number of critical points of $-\psi(x)$ on $\R^n\setminus \Acal$ with index $i$.

Let $n_{\C\setminus \R}$ be the number of non-real critical points of $-\psi(x)$ and $n_{\text{odd}}$ (resp. $n_{\text{even}}$) be the number of real critical points of $-\psi(x)$ with odd (resp. even) index.
Then,
$$
n_{\C\setminus \R}+ n_{\text{odd}}+n_{\text{even}} =|\chi\big(\C^n \setminus (\Acal\cup V(g))\big)| = \sum_{i=0}^n\mu^i(f^hx_0),
$$
and 
$$
n_{\text{even}}-n_{\text{odd}} = \chi(\R^n \setminus \Acal).
$$
Each region of $\R^n\setminus\Acal$ contains at least one local minimum of $-\psi(x)$ and a local minimum has index 0. 
Therefore, 
$$
\# \text{ regions of } \R^n \setminus \Acal \leq n_{\text{even}} \leq \frac{1}{2}\bigl( |\chi\big(\C^n \setminus(\Acal\cup V(g))\big)| + \chi(\R^n \setminus \Acal) \bigr).
$$
\end{proof}
\begin{remark}
From the above proof, we also see that $\# \text{ regions of } \R^n \setminus \Acal$ is bounded above by $\frac{1}{2}(n_\R+\chi(\R^n\setminus \Acal))$ where $n_\R$ is the number of real critical points for $\psi(x).$
Moreover, the equalities in Theorem \ref{thm: upper bound for number of regions using real euler char} hold when all critical points of $\psi(x)$ are real and there is exactly one even index critical point per region.
\end{remark}

We conclude this section with two examples regarding hypersurface arrangements in $\R^2$.

\begin{example}(Generic plane curve of degree $d$)

Let $f_1\in \R[x_1,x_2]$ be a generic real plane curve of degree $d$, $f_2$ be a generic real line and $g\in \R[x_1,x_2]$ be a generic real plane curve of degree 2 such that $g(x)>0$ for all $x\in \R^2$.
The plane curve $V(f_1)$ consists of ovals if $d$ is even and ovals with one pseudoline if $d$ is odd. 
By genericity, $V(f_1)$ is nonsingular and has no self-intersection points.
Each oval has Euler characteristic 0 and a pseudoline has Euler characteristic 1.
A line intersects $V(f_1)$ in at most $d$ real points.
We thus have $\chi(\R^2 \setminus V(f_1f_2)) \leq 1-0-(-1)+ d = d+2$ when $d$ is even and $\chi(\R^2\setminus V(f_1)) =1-1-(-1) +d = d+1$ when $d$ is odd.
For $|\chi(\C^n \setminus (V(f_1)\cup V(f_2) \cup  V(g))|$, it is equal to the coefficient of $z^2$ in $\frac{(1-z)^2}{(1+z)(1+(d-1)z)}$ by Lemma \ref{lem: bn}, which equals $d^2+d+2$.
Hence, 
$$
\# \text{ regions of } \R^2 \setminus V(f_1f_2) \leq \lceil \frac{d^2+d+2 +d+ 1}{2} \rceil = \lceil \frac{d^2+3d+3}{2} \rceil.
$$

% On the other hand, by Harnack's Curve Theorem \cite{harnack1876ueber}, the number of connected components for $V(f_1)$ is bounded above by $\frac{(d-1)(d-2)}{2}+1$ and since the number of regions for $\R^2\setminus V(f_1)$ is the number of connected components for $V(f_1)$ plus one and a line divides each region into at most two regions, we obtain 
% $$
% \# \text{ regions of } \R^2 \setminus V(f) \leq \frac{(d-1)(d-2)}{2}+2.
% $$
% We note that the two upper bounds are asymptotically close as $d$ tends to infinity.

\end{example}

\begin{example}\label{ex:concentric circles} (Concentric circles)

Consider $k$ concentric circles defined by $f_1 = x_1^2+x_2^2-r_1^2,\ldots,f_k = x_1^2+x_2^2-r_k^2$ with distinct $r_1,\ldots, r_k \in \R_{>0}$. 
Define $f_{k+1} = x_1$. 
Denote $\Acal = \cup_{i=1}^{k+1} V(f_i)$, then the number of regions for $\R^2 \setminus
\Acal$ is $2k+2$.
Each region is contractible so $\chi(\R^2\setminus \Acal)=2k+2$.

Let $g\in \R[x_1,x_2]$ be a generic real plane curve of degree 2 such that $g(x)>0$ for all $x\in \R^2$.
We compute $|\chi\big(\C^2\setminus(\Acal\cup V(g))\big)| = (-1)^2 \chi(\C^2\setminus(\Acal\cup V(g)))$.
For any irreducible degree two homogeneous polynomial $p(x_0,x_1,x_2)\in \C[x_0,x_1,x_2]$, $\chi(V(p))$ equals the coefficient of $z^2$ in $\frac{2z}{(1-z)^2(1+z)}$ which is 2, by Lemma \ref{lem: hirzeburch}. 
We thus obtain $\chi(V(f_i))=\chi(V(g))=0$ for $i=1,\ldots,k$ since the homogenization of each $f_i$ and $g$ is irreducible of degree two and the corresponding hypersurface in $\CP^2$ intersects the hyperplane at infinity at two distinct points. 
By inclusion-exclusion for the Euler characteristic, we have 
\begin{align*}
\chi(\Acal\cup V(g))  &= \sum_{i=1}^{k+1}\chi(V(f_i))+\chi(V(g))- \sum_{i=1}^{k+1} \chi(V(f_i)\cap V(g)) - \sum_{i=1}^k \chi(V(f_i) \cap V(f_{k+1}) )\\ &= 1-4k-2-2k = -6k-1.
\end{align*}
Therefore, $|\chi\big(\C^2\setminus(\Acal\cup V(g))\big)|= 6k+2$ and the inequality of Theorem \ref{thm: upper bound for number of regions using real euler char} gives
$$
\# \text{ regions of } \R^2 \setminus \cup_{i=1}^{k+1} V(f_i) \leq  \frac{2k+2+6k+2}{2}  = 4k+2.
$$
The bound we obtain is not tight because not all complex critical points for $\psi(x)$ are real. 
For example, when $k=1, f_1(x)= x_1^2+x_2^2-1, f_2(x) = x_1$ and $g(x) = (x_1-1)^2 + (x_2-2)^2 + (x_1-x_2)^2 +1$, the function $\psi(x) = 2\log|g(x)|- \log|f_1(x)|- \log|f_2(x)|$ has 8 critical points but only 6 of them are real.
It is worth noting that the first inequality in Theorem \ref{thm: number of regions bounded by complex euler} yields an asymptotically tighter bound of 
$4k+2$ compared to the second inequality of the same theorem, which gives an upper bound of $(2k+1)^2+(2k+1)+1 = 4k^2+6k+3$.
\end{example}

\section*{Acknowledgments}
We would like to Lukas Kühne, Chiara Meroni, Leonie Mühlherr, Bernd Sturmfels, Simon Telen for helpful discussions. 
We would also like to note that Simon Telen independently provided an alternative proof of Theorem \ref{thm: one real point per region}.
The second author would like to thank Max Planck Institute for Mathematics in the Sciences for its hospitality and for providing the inspiring environment where the collaboration of this paper started.

\bibliographystyle{alpha}
\bibliography{biblio}
\footnotesize{
\noindent \textsc{Bernhard Reinke}\\
\textsc{Max Planck Institute for Mathematics in the Sciences\\ 
Inselstraße 22, Leipzig 04103} \\
\url{bernhard.reinke@mis.mpg.de} \\

\noindent \textsc{Kexin Wang} \\
\textsc{ Harvard University \\
29 Oxford Street, Cambridge, MA 02138} \\
\url{kexin_wang@g.harvard.edu} \\
}
\end{document}